\documentclass[preprint,12pt]{jams-l}

\usepackage{setspace}
\usepackage{indentfirst}

\usepackage{amsmath,amsfonts,theorem,epsfig,amssymb,overcite}
\usepackage{subfigure}
\usepackage{float}
\usepackage{cite}

\makeatletter
        \renewcommand{\@biblabel}[1]{\textsuperscript{#1}}
        \makeatother

\begin{document}

\title{Efficient implementation of the hybrid method for stochastic simulation of biochemical systems}

\author{Shuo Wang}
\address{Departments of Computer Science, Virginia Tech, Blacksburg, VA 24060}
\email{wangshuo@cs.vt.edu}

\author{Yang Pu}
\address{Departments of Computer Science, Virginia Tech, Blacksburg, VA 24060}
\email{ypu@cs.vt.edu}

\author{Layne Watson}
\address{Departments of Computer Science, Virginia Tech, Blacksburg, VA 24060}
\email{lwatson@cs.vt.edu}

\author{Yang Cao}
\address{Departments of Computer Science, Virginia Tech, Blacksburg, VA 24060}
\email{ycao@cs.vt.edu}

\date{}



\begin{abstract}
Stochastic effect in cellular systems has been an important topic in systems biology. Stochastic
modeling and simulation methods are important tools to study stochastic effect. Given the low efficiency of stochastic simulation algorithms,
the hybrid method, which combines an ordinary differential equation (ODE) system with a stochastic
chemically reacting system, shows its unique advantages in the modeling and simulation of biochemical systems.
The efficiency of hybrid method is usually limited by reactions in the stochastic subsystem,
which are modeled and simulated using Gillespie's framework and frequently interrupt the integration of the ODE subsystem.
In this paper we develop an efficient implementation approach for the hybrid method coupled with traditional ODE solvers. We also 
compare the efficiency of hybrid methods with three widely used ODE solvers RADAU5, DASSL, and DLSODAR.
Numerical experiments  with three biochemical models are presented.
A detailed discussion is presented for the performances of three ODE solvers.
\end{abstract}

\maketitle

\section{Introduction} \label{sec:intro}

Living cells are regulated by complicated biological networks. To study the mechanism of these networks,
systems biologists often use chemically reacting systems and ordinary differential equations (ODEs)
to model the behavior of cellular networks. In living cells, however, reacting molecules participate in chemical reactions
in a discrete and stochastic manner rather than a continuous and deterministic one.
This stochastic feature results in noise, which sometimes has significant impacts on living cells.
In order to take the noise into consideration, stochastic simulation methods are needed. Among them the most famous one is Gillespie's stochastic simulation algorithm (SSA)~\cite{ssa1, ssa2} for well-stirred chemically reacting systems. 
However, as models become more complex, the SSA often becomes computationally expensive because it simulates every reaction firing. In biochemical systems, some reactions involve species of large populations or have large reaction rates. These reactions tend to dominate the CPU time in an SSA simulation.
In order to speed up the SSA, several approximation strategies have been proposed. One group of approximation methods tries to take advantage of the multiscale property among reactions: Some reactions are much faster than others and these fast reactions may reach partial equilibrium~\cite{ssSSA},
or quasi steady state~\cite{QSSA} during each slow reaction firing. These methods include the quasi steady state SSA method~\cite{QSSA}
and the slow-scale SSA method~\cite{ssSSA}.
Another group of approximation algorithms was proposed to take advantage of the multiscale characteristics among species  populations.
The tau-leaping methods~\cite{Gillespie2001} belong to this group. To combine the tau-leaping method and the SSA,
Cao et al.~\cite{caoTauLeaping} proposed a hybrid method to partition species into two sets according to their population numbers. Reactions involving species whose population are less than a threshold are simulated by the SSA and other reactions are simulated by the tau-leaping method.

Of course, multiscale features among reactions and species populations are connected. Often fast reactions are related to species with large populations. Based on this observation, Haseltine and Rawlings~\cite{SlowFast} developed a strategy to partition a system into a slow reaction group and a fast reaction group.
The fast reaction group consists of fast reactions whose reactants are of large populations and
are modeled by ODEs or chemical Langevin equations (CLEs).
The rest reactions, either the reaction is slow or one of the involved species is in a small population,
are put in the slow reaction group and are modeled and simulated by the SSA.
Recently a new partitioning strategy was proposed to put all fast reactions, no matter its reactants are in large or small populations, into the fast reaction group and simulated by ODEs.
This strategy made Haseltine and Rawling's hybrid method more efficient and practical when applied to biochemical systems. On the other hand,
the implementation details for the hybrid method were not discussed in literature at all although it is critical for the success of the hybrid method.

To implement Haseltine and Rawling's hybrid method, numerical solvers for stiff ODEs are always needed due to the multiscale features.
There are two major groups of ODE solvers for stiff problems: one step methods and multistep methods. The former group include the famous Runge-Kutta methods, while the latter group include Adams methods and BDF methods, which many popular stiff ODE solvers are based on.
ODE solvers based on multistep methods have received more attention in literature since these methods usually lead to smaller linear system size
than the one step methods and thus are considered faster than one step methods. However, when applying the hybrid method to
biochemical systems, there is a very special feature. ODE integration is often interrupted by discrete slow reactions. In that case,
ODE solvers based on multistep methods have to reduce orders and stepsizes to handle the discontinuity resulted from discrete slow reactions.
When slow reactions fire frequently, it is a difficult job for multistep methods to maintain its high order and large stepsizes. Then
the efficiency of multistep methods will be heavily affected.
On the contrary, step sizes of one step methods will not be affected by the discontinuity resulted from discrete slow reactions.
They are usually greater than stepsizes selected by multistep methods. Therefore, we conjecture that one step methods should be a better choice when the
hybrid method is concerned. But before we draw a solid conclusion, we need to study what is the most efficient way to implement the hybrid
method using one step ODE solvers, and we need to perform numerical experiments on typical biological models to verify/reject our conjecture.

In this paper, we present our study of the implementation details related to the hybrid method.
We first developed an efficient strategy for the root finding function in any ODE solver. This strategy utilizes a special feature in the hybrid method: the event function is monotonously increasing. With this strategy, the event handling will be more efficient and can be applied to both one step methods and multistep methods.
To investigate the accuracy and efficiency of hybrid method in terms of different ODE solvers, we did numerical experiments based on three biological models.
The ODE solvers we tested are RADAU5~\cite{SolvingII}, DASSL~\cite{DASSL}, and DLSODAR~\cite{DLintro},
where RADAU5 is an implicit one step Runge-Kutta method ODE solver, DASSL and DLSODAR are ODE solvers based on multistep methods.
Numerical results demonstrate that the hybrid method with these three ODE solves have similar accuracy but different efficiency performance.
For the two multistep ODE solvers, DLSODAR has better performance. 
In the first two experiments, RADAU5 shows its advantage when the integration process is interrupted frequently by slow reactions. In the third experiments, RADAU5 does not perform as well as expected when slow reaction firing frequency increases. Performance data collected in numerical experiments indicate that the drastically growing CPU time on Jacobian matrix calculation largely weakens RADAU5's efficiency.

This paper is organized as follows. In Section \ref{sec:background}, we briefly review Gillespie's SSA and introduce the hybrid method and ODE solvers applied in our numerical experiments. In Section \ref{sec:inverse}, we propose an efficient event handling function in terms of inverse interpolation. The accuracy of this method is also discussed. Three numerical experiments and corresponding results are shown in Section \ref{sec:example}. In Section \ref{sec:discussion}, we give our conclusion and further discussion.

\section{Background}\label{sec:background}
\subsection{Stochastic simulation algorithm}

We consider a well-stirred system of $N$ chemical species $\{S_1,\dots,S_N\}$ and $M$ chemical reactions $\{R_1,\dots,R_M\}$. We define the current state of system as a vector $\textbf{x}=(x_1,\dots,x_N)$, where $x_i$ is the number of $S_i$ molecules. Each reaction channel $R_j$ is characterized by its propensity function $a_j(\textbf{x})$ and its state-change vector $\textbf{v}_j=(v_{1j},\dots,v_{Nj})$. The probability of reaction $R_j$ firing in the next infinitesimal time $dt$ is given as $a_j(\textbf{x})dt$, and $v_{ij}$  $S_i$ molecules will be changed by one $R_j$ reaction. The procedure of SSA is presented as follows.

\begin{itemize}
\item[1)] In state $\textbf{x}$ at time $t$, compute propensity values $a_j(\textbf(x))$, and the sum $a_0(\textbf{x}) \equiv \sum^M_{j=1} a_j(\textbf{x})$.
\item[2)] Generate a time increment $\tau$ as a sample of the exponential random variable with mean $1/a_0(\textbf{x})$.
\item[3)] Generate a reaction index $j$ based on the probability function $a_j(\textbf{x})/a_0(\textbf{x}) (j=1,\dots,M)$.
\item[4)] Update time $t\leftarrow t+\tau$ and state $\textbf{x}\leftarrow\textbf{x}+\textbf{v}_j$.
\item[5)] Return to step 1) if stopping condition is not reached.
\end{itemize}

Although SSA is mathematically exact, the time cost by SSA may become very expensive when the size of a biochemical system increases, since SSA should simulate every reaction event and the time increment $\tau$ may be very small corresponding to a large $a_0(\textbf{x})$. Thus it is necessary to develop more efficient simulation methods.

\subsection{Hybrid method}

The hybrid method follows the strategy proposed by Haseltine and Rawlings.
Given a system of $N$ species $\{S_1,\dots,S_N\}$ and $M$ reactions $\{R_1,\dots,R_M\}$, these $M$ reactions are partitioned into two subsets. Reactions in subset $S_{fast}$ are formulated by ODEs, and reactions in subset $S_{slow}$ are simulated with SSA. Let $a_i(\textbf{x},t)$ be the propensity function of $i$-th reaction in $S_{slow}$, where $\textbf{x}$ is the state at time $t$, and its state-change vector $\textbf{v}_i=\{v_{i1},\dots,v_{iN}\}$. Let $\tau$ be the jump interval of the next stochastic reaction, and $\mu$ be its reaction index. Set $t=0$. The hybrid method simulate the system as follows:

\begin{itemize}
\item[1)] Generate two uniform random numbers $r_1$ and $r_2$ in $U(0,1)$.
\item[2)] Integrate the ODE system and the integral equation:
\begin{equation}
\int^{t+\tau}_{t}a_{tot}(\textbf{x},s)ds+\log(r_1) = 0,
\label{eq:SolveTau}
\end{equation}
where $a_{tot}(\textbf{x},t)$ is the sum of propensities of reactions in $S_{slow}$.
\item[3)] Determine $\mu$ as the smallest integer satisfying
\begin{equation}
\sum^{\mu}_{i=1}a_i(\textbf{x},t)>r_2a_{tot}(\textbf{x},t).
\label{eq:SloveMu}
\end{equation}
\item[4)] Update $\textbf{x}\leftarrow \textbf{x}+\textbf{v}_{\mu}$.
\item[5)] Return to step 1) if stopping condition is not reached.
\end{itemize}

In step 2), different from Haseltine and Rawlings's strategy of adding a propensity of ``no reaction'', we choose a more direct strategy to get the jump interval $\tau$. Suppose that the ODE system is given by
\begin{equation}
\textbf{x}'=f(\textbf{x}).
\label{eq:ODEsys}
\end{equation}
We simply add an integration variable $z$ and an equation
\begin{equation}
z'=a_{tot}(\textbf{x}), z(t)=\log(r_1), 
\label{eq:Zsys}
\end{equation}
where we note that $\log(r_1)$ is negative. 
During the simulation, each step starts at time $t$ and numerically integrates ODEs (\ref{eq:ODEsys}) and (\ref{eq:Zsys}). When $z(t+\tau)=0$, the ODE integration stops, where $\tau$ is the solution to (\ref{eq:SolveTau}). This integration process can be conveniently handled by standard ODE solvers combined with root-finding functions. Note that since $z$ is an integration variable, one may choose to omit it from the error control mechanism~\cite{DASSL}. Adding this extra variable will not greatly affect the efficiency of ODE solution. 

\subsection{ODE solvers}

When applying the hybrid method to simulate biochemical systems, all reactions in the subset $S_{fast}$ and the equation (\ref{eq:Zsys}) need to be solved by an ODE solver. Because of the multiscale property of biochemical systems, ODE solvers for stiff problems are usually preferred. There are two types of widely-used ODE solvers (in Fortran) for stiff problems, one based on one step methods and the other based on multistep methods.

For one step methods, in each step to integrate from $t_n$ to $t_{n+1}$, the integration only uses the value $\textbf{x}(t_n)$.
For stiff systems, implicit methods are usually used.
RADAU IIA (RADAU5) is an implicit Runge-Kutta method of stage three and order five for stiff ODEs, which was originally
proposed by Axelsson in 1969~\cite{SolvingII}. In each step, RADAU5 uses the current state at time $t$ to calculate the state at time $t+h$, where $h$ is the step size. The Fortran code we used for RADAU5 was written by E. Hairer and G. Wanner~\cite{SolvingII}. It offers efficient step size control strategy and continuous output, but no root finding function. However, this work can be done with the assistant of function CONTR5 when calling function SOLOUT. More details can be found in the user guide of the RADAU5 code.

DASSL and DLSODAR are widely-used multistep ODE solvers for stiff problems. In contrast to one step methods, multistep methods use state values at several previous time points to calculate the state value at $t_{n+1}$. Adams method and BDF (backward differentiation formulas) method are two popular multistep methods.
The idea of Adams method is to approximate $\textbf{x}(t_{n+1})$ using derivatives from several previous time points, while the BDF method
is based on an interpolation polynomial using state values at previous time points to approximate the derivative at $t_{n+1}$. 
DASSL was developed by Linda R. Petzold \cite{DASSL}. It follows the idea of BDF methods. Similar as RADAU5, DASSL also does not directly
support root handling function.
DLSODAR is an ODE solver developed by Alan C. Hindmarsh and Linda R. Petzold\cite{DLintro}. A great feature of DLSODAR is 
 that it can automatically switch between Adams methods and BDF methods based on characteristics of underlying ODEs. This property allows DLSODAR to have a better performance compared with DASSL\cite{DLAutoSelect}. It also offers root-finding function. Thus it is more convenient for event handling.

\section{Inverse interpolation}\label{sec:inverse}

In the hybrid method, the event that a slow reaction in $S_{slow}$ fires is important and has significant impact on the accuracy of the simulations.
Assume the integration for the equation (\ref{eq:Zsys}) starts at time $t$. For a slow reaction to fire at $t_{event}$, we should have
$z(t_{event})=0$.
In order to obtain the solution $t_{event}$, ODE solvers have to equip with root-finding functions. DLSODAR has a root-finder based on Illinois Algorithm\cite{DLrootfind}, a variant of secant method, which returns the leftmost solution in terms of iterations. But for RADAU5 and DASSL, root finding functions
are not directly provided.
In order to compare the performance of RADAU5, DASSL, and DLSODAR for the hybrid method, we need to develop an efficient root finding mechanism
for them.
Generally, a root finding function in an ODE solver needs to develop an
interpolation polynomial using state values at time points (for multistep methods) or stage points (for Runge-Kutta methods). However, root finding function in the hybrid method has a very special feature, which makes the inverse interpolation method a very good choice.
In equation (\ref{eq:Zsys}), the derivative of $z(t)$ is given by $a_{tot}(\textbf{x})$,
the sum of propensities of slow reactions. Thus $z'(t)$ is always positive. It implies that $z(t)$ is a monotonous function with respect to time $t$.
So there must exist an inverse function for $z(t)$. We denote the inverse by $t = T(z)$.
Thus $t_{event}$, the solution  of $z(t) = 0$,  can be given by $t_{event} = T(0)$. In the following we give the details on how the inverse interpolation method is implemented for RADAU5 and DASSL.

\subsection{Inverse interpolation in RADAU5}

RADAU AII is a three-stage five-order implicit Runge-Kutta method. For an ordinary differential equation
\begin{equation}
y'(t)=f(t,y),
\label{eq:ysys}
\end{equation}
the value at time $t+h$ is computed based on the calculated value at time $t$ in the previous step and corresponding formulas are as follows
\begin{equation}
\begin{array}{rclc}
g_i&=&y(t)+h\displaystyle{\sum_{j=1}^{3}} a_{ij}f(t+c_{j}h,g_j),& i=1,2,3,\\
y(t+h)&=&y(t)+h\displaystyle{\sum^3_{j=1}}b_{j}f(t+c_{j}h,g_j),&
\end{array}
\label{eq:radau}
\end{equation}
where $h$ is the step size, $a_{ij}$, $b_j$ and $c_j$ are coefficients in RADAU AII\cite{SolvingII}. For $i=1,2$, $g_i$ is an approximation of $y(t+c_{j}h)$. For $i=3$, $a_{3j}$ shares the same value with $b_j$ and $c_3$ is equal to $1$, therefore $g_3$ gives the value of $y$ at time $t+h$. All $g_i$'s are calculated using implicit method.

When running the hybrid method with RADAU5, the integration process showed in equation (\ref{eq:radau}) is applied to both fast reactions described as ODEs and the equation (\ref{eq:Zsys}) step by step when $z$ is less than $0$. If $z\geq 0$ is detected at time $t+h$, we call a root-finding function using inverse interpolation to find the time $t < t_{event}\leq t+h$ by solving $z(t_{event})=0$. 
For example, the inverse interpolation function $T_1(z)$ can be constructed in the Lagrange form as
\begin{equation}
\begin{array}{ll}
T_1(z)=& \frac{(z-g_1)(z-g_2)(z-g_3)}{(z_t-g_1)(z_t-g_2)(z_t-g_3)}t \\
	& +  \frac{(z-z_t)(z-g_2)(z-g_3)}{(g_1-z_t)(g_1-g_2)(g_1-g_3)}(t+c_1h)\\
       & + \frac{(z-z_t)(z-g_1)(z-g_3)}{(g_2-z_t)(g_2-g_1)(g_2-g_3)}(t+c_2h)\\
       & +\frac{(z-z_t)(z-g_1)(z-g_2)}{(g_3-z_t)(g_3-g_1)(g_3-g_2)}(t+h)
\end{array}
\label{eq:radau_inverse}
\end{equation}
The solution to $z(t_{event})=0$ can be obtained as $t_{event}=T_1(0)$.

\subsection{Inverse interpolation in DASSL}

DASSL is an ODE solve based on the BDF method. The solution to equation (\ref{eq:ysys}) at step $n$ can be obtained as
\begin{equation}
y_n=\displaystyle{\sum_{i=1}^q}\alpha_iy_{n-i}+h\beta_0f(t_n,y_n)
\label{eq:BDF}
\end{equation}
where $\alpha_i$ and $\beta_0$ are coefficients, $h$ is the step size, and $q$ is the order indicating how many previous time points are needed. In DASSAL, the order $q$ varies from one to five.

Different from the one step method used in RADAU5, DASSL does not calculate the value between $y_{n-1}$ and $y_{n}$. Therefore we have to adopt a little different inverse interpolation strategy in DASSL. During the process of integration, we use two vectors of length six to record the values of $z$ for equation (\ref{eq:Zsys}) and their corresponding time $t$. When $z_n>0$ is detected, we apply $z_n$ and previous $q$ results $z_{n-1},\dots,z_{n-q}$ as independent variables and $t_n,\dots,t_{n-q}$ as dependent variables to construct a polynomial $T_2(z)$ with highest order $q$. The solution $t_{event}$ where $z(t_{event})=0$ is got at $t_{event}=T_2(0)$.

\subsection{Accuracy of inverse interpolation}

In this part, we will briefly discuss the accuracy of inverse interpolation in both RADAU5 and DASSL.

First, we give an accuracy analysis of inverse interpolation for RADAU5. Let $z(t)$ be the exact value of equation (\ref{eq:Zsys}) at time $t$ and $z_t$ is the numerical solution to $z(t)$. According to Lemma 7.5 in \cite{SolvingI}, we have
\begin{equation}
\begin{array}{rcl}
 z_t &=& z(t)+O(h^5),\\
 g_1 &=& z(t+c_1h)+O(h^3),\\
 g_2 &=& z(t+c_2h)+O(h^3),\\
 g_3 &=& z(t+h)+O(h^5).
\end{array}
\label{eq:gOrder}
\end{equation}
In order to simplify this analysis, we assume that $z_t<g_1<g_2<g_3$. Let $t^*$ be the exact time increment where $z(t^*)=0$ and $\hat{t}$ be the approximation obtained from equation (\ref{eq:radau_inverse}). If we replace $z_t$ and $g_i$ by $z(t)$ and $z(t+c_ih)$ in $T_1(0)$ respectively, we have an approximation denoted by $t_z$. The error can be represented as
\begin{equation}
err(\hat{t})=(\hat{t}-t_z) + (t_z -t^*)
\label{eq:tauErr}
\end{equation}
This error is composed of two parts. The first part on the right hand side gives the error from approximating $z(t)$ and $z(t+c_ih)$ by $z_t$ and $g_i$. The second part is the interpolation error.
To obtain the error bound for the first part, we use the Lagrange form \eqref{eq:radau_inverse} of $T_1(z)$ and let $z = 0$.
Let $z_0$ denote $z(t)$ and $z_i$ denote $z(t+c_ih)$, by ignoring the higher order of $h$, we have
\begin{equation}
\begin{array}{rcl}
\frac{g_1g_2g_3}{(z_t-g_1)(z_t-g_2)(z_t-g_3)}t &=& \frac{(z_1+O(h^3))(z_2+O(h^3))(z_3+O(h^5))}{\displaystyle{\prod^3_{i=1}}(z_0-z_i+O(h^3))}t\\
                                      &=& \frac{z_1z_2z_3}{(z_0-z_1)(z_0-z_2)(z_0-z_3)}t+O(h^3)
\end{array}
\label{eq:gErr}
\end{equation}
For the rest two items in $T_1(0)$, we have similar results as in $\eqref{eq:gErr}$. Therefore, $\hat{t}-t_z=O(h^3)$. Moreover, the inverse interpolation based on $z(t+c_ih)$ has an error $O(h^4)$ that indicates $t_z -t^*=O(h^4)$. We can get the total error $err(\hat{t})=O(h^3)$.

We can apply a similar accuracy analysis to the inverse interpolation in DASSL. We define $z(t_i)$ as the exact solution to equation (\ref{eq:Zsys}) at time $t_i$ and $z_i$ as the numerical solution calculated from equation (\ref{eq:BDF}) using $q_i$ previous solutions. Let $t^*$ be the absolute time where $z(t^*)=0$, $\hat{t}$ be the approximation given by $T_2(0)$ and $t_z$ be the solution after replacing $z_i$ by $z(t_i)$ in $T_2$. Theory 2.4 in \cite{SolvingI}gives that $z_i=z(t_i)+O(h^{q_i})$. We thus have $\hat{t}-t_z=O(h^{min(q_i)})$. On the other hand, we assume that $q+1$ exact values are used to construct inverse interpolation polynomial. We have $t_z-t^*=O(h^q)$. Therefore the order of solution $\hat{t}$ is $min(q_i,q)$. In DASSL, $q$ varies from one to five. The inverse interpolation in DASSL has a maximum order five.

\section{Numerical examples}\label{sec:example}

In a hybrid method simulation, each time when there is a slow reaction firing, the ODE solver should stop and restart after the state is updated. For a multistep ODE solver like DASSL, the restart procedure may have great overhead if slow reactions fire too often because it has to reduce to a first order method with a small stepsize. This process has a heavy impact on the efficiency.  For a Runge-Kutta solver like RADAU5, the impact should be small.
In this section we compare the accuracy and efficiency of hybrid methods with three different ODE solvers: RADAU5, DASSL, and DLSODAR.
Numerical results for three examples are present. One example is a simple model involving three species and three reactions.
The second example is a gene regulation process, and the third one is a cell cycle model presented in Liu et al.\cite{Zhen}.
Because they are stochastic simulations, the number of species and simulation CPU time vary in each run.
We therefore recorded the number of slow reactions firing times and CPU time for each simulation.
Accuracy and efficiency are measured in terms of average results, respectively. We set absolute error tolerance as $10^{-6}$ and relative error tolerance as $10^{-3}$ for all three ODE solvers. The simulations were performed on a 2.30GHz Intel Core Liinux workstation.

\subsection{Example one}

\begin{figure}
  \centering
  \subfigure[Average firing times of the third reaction on example one]{
    \label{fig:toy_num}
    \epsfig{file=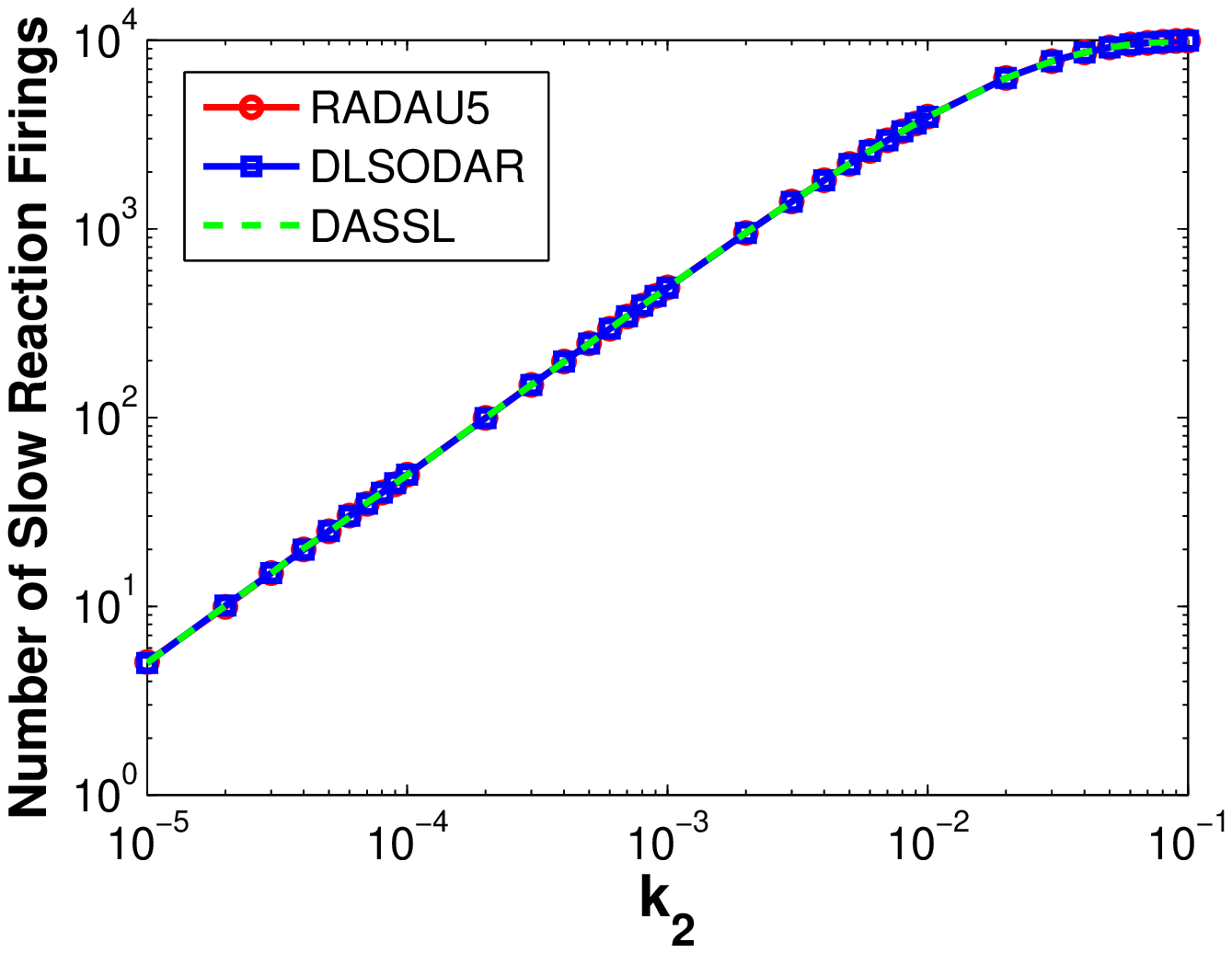, width=2.5in}}
  \hspace{.01in}
  \subfigure[Average CPU time of hybrid method with three ODE solvers one example one]{
    \label{fig:toy_time}
    \epsfig{file=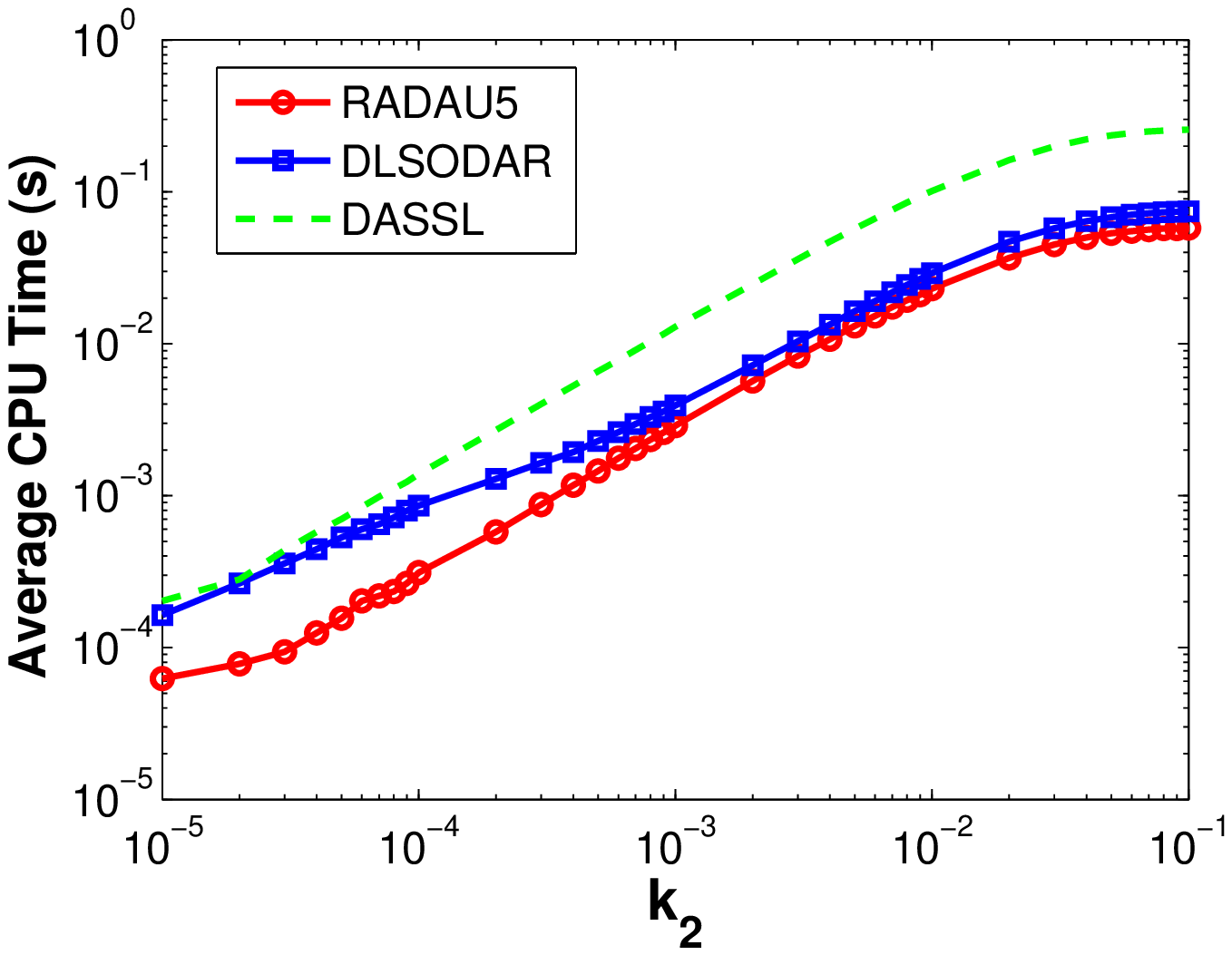, width=2.5in}}
  \caption{Accuracy and efficiency comparison of hybrid method with three ODE solvers on example one}
  \label{fig:toy}
\end{figure}

We first consider a simple system
\begin{equation}
\begin{array}{ccc}
S_1 &\xrightarrow{k_1} & S_2, \\
S_1 &\xleftarrow{k_{-1}} & S_2, \\
S_2 &\xrightarrow{k_2} & S_3,
\end{array}
\label{eq:toy}
\end{equation}
where $k_1=0.5$, $k_{-1}=1.5$ and $k_2$ varies from $10^{-5}$ to $10^{-1}$. The initial values of three species are $x_1(0)=7500$, $x_{2}(0)=2500$, and $x_{3}(0)=0$. Here the first two reversible reactions are considered as fast reactions and simulated by ODE solvers, and the last reaction is a slow reaction and handled by SSA. This system were simulated 1,000 times with $T=200$ as the final time using ODE solvers RADAU5, DASSL, and DLSODAR, respectively. The average CPU time and the number of  firing times for the third reaction are shown in Figure \ref{fig:toy}.

Figure \ref{fig:toy_num} shows the mean value of the third reaction firing number and Figure \ref{fig:toy_time} gives the mean value of CPU times, corresponding to different $k_2$ values. In Figure \ref{fig:toy_num}, the three curves match with each other very well. It suggests that given the same error tolerance, hybrid methods with RADAU5, DASSL and DLSODAR have similar accuracy for this model. But from Figure \ref{fig:toy_time}, it is easy to see that the average running times of the hybrid method with the three ODE solvers are different.
RADAU5 performs the best for this example, while DASSL is the worst. For DLSODAR, when $k_2<10^{-4}$ the curve is close to the one for DASSL and when $k_2>10^{-3}$ the curve is close to the one for RADAU5. But there is no crossover. So for this simple system, RADAU5 is the best choice for the hybrid method.

\subsection{Gene regulation model}

\begin{figure}
 \centering
 \epsfig{width=2in, file=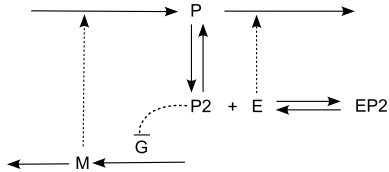}
 \caption{Diagram for gene regulation model}
 \label{fig:gene_model}
\end{figure}

The second model is a gene regulation model, which was first given in the appendix of paper \cite{Zhen}. The diagram for this model is shown in Figure \ref{fig:gene_model}.

In this model, the synthesis of protein $p$ is controlled by the expression of gene $g$. Protein $p$ also can form homodimer $p2$, which can  bind to the promotor site of its own gene and repress the expression (a negative feedback). In addition, $p2$ also binds and inhibits the enzyme $E$, which catalyzes the degradation of protein $p$ (a positive feedback). Oscillation of total level of protein ($p_{total}=p+2*p2$) occurs during the simulation.
Reactions, rate constants and partition are given in Table \ref{tab:gene_reaction}. The initial value of species are shown in Table \ref{tab:gene_init}.

\begin{table}
    \centering
    \caption{Reactions in the Gene Regulation Model}
    \begin{tabular}{lll}
        \hline
        Reaction   &            Rate Constant    &     Partition\\
        \hline
        $g\rightarrow g+m$      &   $k_1=0.5$        &      SSA\\
        $m\rightarrow m+p$      &   $k_2=60$         &      ODE\\
        $p\rightarrow \phi$     &   $k_3=0.05$       &      ODE\\
        $m\rightarrow \phi$     &   $k_4=0.05$       &      SSA\\
        $p+p\rightarrow p2$     &   $k_5=0.001$      &      ODE\\
        $p2\rightarrow p+p$     &   $k_6=40$         &      ODE\\
        $p2+g\rightarrow gi$    &   $k_7=40$         &      ODE\\
        $gi\rightarrow p2+g$    &   $k_8=1000$       &      ODE\\
        $E+p2\rightarrow Ep2$   &   $k_9=0.4$        &      ODE\\
        $Ep2\rightarrow E+p2$   &   $k_{10}=0.5$     &      ODE\\
        $E+p\rightarrow Ep$     &   $k_{11}=0.2$     &      ODE\\
        $Ep\rightarrow E+p$     &   $k_{12}=10$      &      ODE\\
        $Ep\rightarrow E$       &   $k_{13}=10$      &      ODE\\
        \hline
    \end{tabular}
    \label{tab:gene_reaction}
\end{table}

\begin{table}
    \centering
    \caption{Initial Value of Reactants in the Gene Regulation Model}
    \begin{tabular}{lcccccccc}
        \hline
        Reactant      & g & gi & m &  p & p2 & E & Ep & Ep2\\
        \hline
        Initial Value & 1 & 0  & 10&1000& 100&100& 0 & 0 \\
        \hline
    \end{tabular}
    \label{tab:gene_init}
\end{table}

\begin{figure}
  \centering
  \subfigure[Average firing times of $m$ synthesis reaction on gene regulation model]{
    \label{fig:gene_syn}
    \epsfig{file=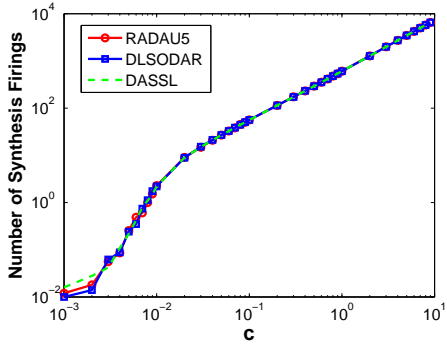, width=2.5in}}
  \hspace{.01in}
  \subfigure[Average firing times of $m$ degradation reaction on gene regulation model]{
    \label{fig:gene_deg}
    \epsfig{file=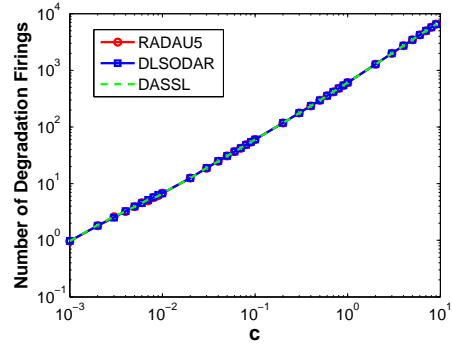, width=2.5in}}
  \caption{Accuracy comparison of hybrid method based on three ODE solvers on gene regulation model}
  \label{fig:gene_acc}
\end{figure}

\begin{figure}
 \centering
 \epsfig{width=3.5in, file=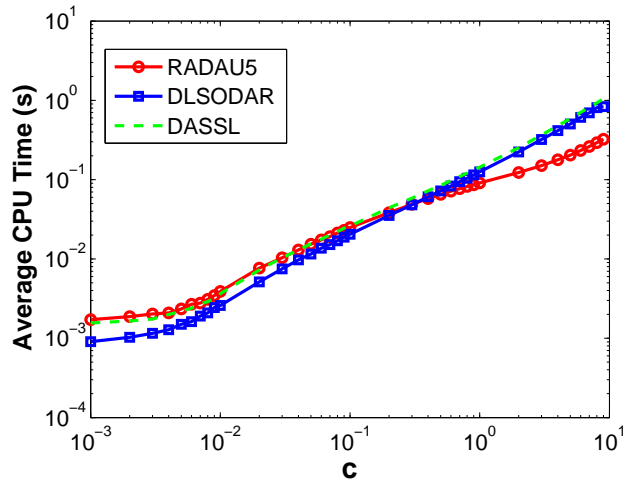}
 \caption{Efficiency comparison of hybrid method based on three ODE solvers on gene regulation model}
 \label{fig:gene_eff}
\end{figure}

In this model, synthesis and degradation of the mRNA $m$ are simulated by the SSA. To test the influence of slow reaction firing frequency on the performance of hybrid methods with different ODE solvers, we introduce a scale factor $c$ and multiply it to $k_1$ and $k_4$. We let $c$ vary from
$10^{-3}$ to $10$. The larger $c$ is, the more frequently an ODE solver has to restart. With $T=2,000$ as the final time, this system were simulated $1,000$ times using the hybrid method with ODE solvers RADAU5, DASSL, and DLSODAR, respectively.  The mean firing numbers of synthesis and degradation and average CPU time are shown in Figure \ref{fig:gene_acc} and Figure \ref{fig:gene_eff}.

From Figure \ref{fig:gene_syn}, we observe that when $c$ is less than $10^{-2}$, the system's stochastic property has an obvious impact on the synthesis of $m$. As $c$ increases, the curves become smooth and match with each other very well. Also in Figure \ref{fig:gene_deg}, we can see that
hybrid methods with the three ODE solvers give similar results. Same as the result shown in Figure \ref{fig:toy_time}, in Figure \ref{fig:gene_eff} the
hybrid method with DLSODAR always uses less time than the one with DASSL. But the hybrid method with RADAU5 has a different performance. When $c$ is less than $0.03$, the average CPU time of RADAU5 is close to the result of DASSL but a little larger. Then the curve for RADAU5 grows slower and has crossovers with the curves for DASSL and DLSODAR at $c=0.03$ and $c=0.3$, respectively. After that, the gap between the red curve and the blue curve becomes larger and larger. It suggests that in this gene regulation model DLSODAR performs the best when less slow reaction fires during the simulation, but
RADAU5 is better when slow reactions fire more frequently.

\subsection{Cell cycle model}

\begin{figure}
 \centering
 \epsfig{width=2in, file=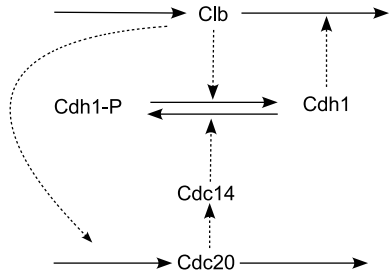}
 \caption{Diagram for cell cycle model}
 \label{fig:cell_cycle}
\end{figure}

The cell cycle is driven by the mutual antagonism between B-type
cyclins (such as Clb2) and G1-stabilizers (such as Cdh1)\cite{cellcycle}.
When B-type cyclins are abundant, they combine with kinase subunits (Cdk1)
to form active protein kinases (e.g., Cdk1-Clb5 and Cdk1-Clb2 in budding yeast)
that promote DNA synthesis and mitosis (S, G2, and M phases of the cell cycle).
When Cdh1 is active, Clb-levels are low, and cells are in the unreplicated phase (G1)
of the DNA replication-division cycle.
The cell cycle control system alternates back-and-forth between
G1 phase (Cdh1 active, Clb-levels low) and the S-G2-M phase (Cdh1 inactive, Clb-levels high).
Tyson and Novak built a three-variable model to describe the bistable switch give by the antagonism between Cdk1-Clb complexes and Cdh1. The diagram is shown in Figure \ref{fig:cell_cycle}.

Let $X$ denote Cdk1-Clb complexes, $Y$ denote Cdh1 and $Z$ denote the lumped reaction of Cdc20 and Cdc14. In this model, $X$ inactivates $Y$ by phosphorylating it, and unphosphorylated $Y$ catalyzes $X$'s degradation. $X$ also promotes the synthesis of $Z$ and $Z$ catalyzes the dephosphorylation
of  phosphorylated $Y$. Therefor, $X$ and $Y$ are involved in a mutual antagonism, which creates a bistable switch.

Based on this three-variable cell cycle model, Liu et al. \cite{Zhen} built a hybrid stochastic cell cycle model by introducing
six stochastic reactions for the synthesis and degradation of three mRNA variables. In this hybrid model, for the ODE part, besides $X$, $Y$, and $Z$, the variable $V$ denotes the volume of cell, and $Y_T$ is the total amount of $Y$ (unphosphorylated and phosphorylated). For the SSA part, $M_x$, $M_y$, and $M_z$ are three mRNA variables for $X$, $Y$, and $Z$, respectively. Details of the model are listed in Table \ref{tab:cell_ODE}, Table \ref{tab:cell_SSA} and Table \ref{tab:cell_para}.

\begin{table}
    \centering
    \caption{ODE system of hybrid cell cycle model}
    \begin{tabular}{lcl}
        \hline
        $\frac{d}{dt}V$ &=& $\mu V$\\
        $\frac{d}{dt}\langle X\rangle$ &=& $k_{sx}M_{x}V-k_{dx}\langle X\rangle-\frac{k_{dxy}\langle X\rangle\langle Y\rangle}{V}$\\
        $\frac{d}{dt}\langle Y_T\rangle$ &=& $k_{sy}M_yV-k_{dy}\langle Y_T\rangle$\\
        $\frac{d}{dt}\langle Y\rangle$ &=& $k_{sy}M_yV-k_{dy}\langle Y\rangle+\frac{(k_{hy}V+k_{hyz}\langle Z\rangle)(\langle Y_T\rangle-\langle Y\rangle)}{J_{hy}V+\langle Y_T\rangle-\langle Y\rangle}$\\
                                       & & $-\frac{k_{pyx}\langle X\rangle\langle Y\rangle}{J_{pyx}V+\langle Y\rangle}$\\
        $\frac{d}{dt}\langle Z\rangle$ &=& $k_{sz}M_zV-k_{dz}\langle Z\rangle$\\
        \hline
    \end{tabular}
    \label{tab:cell_ODE}
\end{table}

\begin{table}
    \centering
    \caption{SSA system of hybrid cell cycle model}
    \begin{tabular}{ll}
        \hline
        Reaction    &    Propensity Function\\
        \hline	
        $\phi\rightarrow M_x$ & $k_{smx}V$\\
        $M_x\rightarrow \phi$ & $k_{dmx}M_x$\\
        $\phi\rightarrow M_y$ & $k_{smy}$\\
        $M_y\rightarrow \phi$ & $k_{dmy}M_y$\\
        $\phi\rightarrow M_z$ & $k_{smz}+\frac{k_{smzx}\langle X\rangle^2}{(J_{smzx}V)^2+\langle X\rangle^2}$\\
        $M_z\rightarrow \phi$ & $k_{dmz}M_z$\\
        \hline
    \end{tabular}
    \label{tab:cell_SSA}
\end{table}

\begin{table}
    \centering
    \caption{Parameter values of hybrid cell cycle model}
    \begin{tabular}{lclc}
        \hline
        Parameter    &    Value ($min^{-1}$)    &     Parameter    &      Value ($fL^{-1}\cdot min^{-1}$)\\
        $\mu$        &        0.006             &     $k_{sx}$     &      1.53\\
        $k_{dx}$     &        0.04              &     $k_{sy}$     &      1.35\\
        $k_{dy}$     &        0.02              &     $k_{hy}$     &      29.7\\
        $k_{hyz}$    &        7.5               &     $k_{sz}$     &      1.35\\
        $k_{pyx}$    &        1.88              &     $k_{smx}$    &      1.04\\
        $k_{dz}$     &        0.1               &     Parameter    &      Value ($fL\cdot min^{-1}$)\\
        $k_{dmx}$    &        3.5               &     $k_{dxy}$    &      0.00741\\
        $k_{smy}$    &        7.0               &     Parameter    &      Value ($fL^{-1}$)\\
        $k_{dmy}$    &        3.5               &     $J_{hy}$     &      5.4\\
        $k_{smz}$    &        0.001             &     $J_{pyx}$    &      5.4\\
        $k_{smzx}$   &        10.0              &     $J_{smzx}$   &      756\\
        $k_{dmz}$    &        0.15              &                  &\\
        \hline
    \end{tabular}
    \label{tab:cell_para}
\end{table}

\begin{figure}
 \centering
 \epsfig{width=3.5in, file=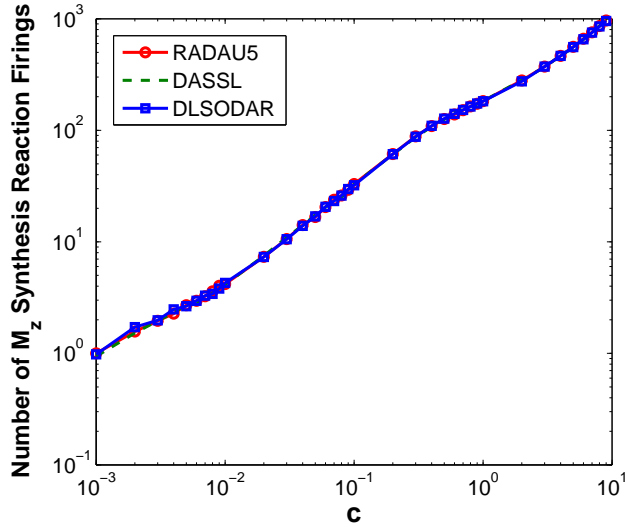}
 \caption{Average firing times of $M_z$ synthesis reaction of hybrid method with three ODE solvers on cell cycle model}
 \label{fig:cell_acc}
\end{figure}

\begin{figure}
 \centering
 \epsfig{width=3.5in, file=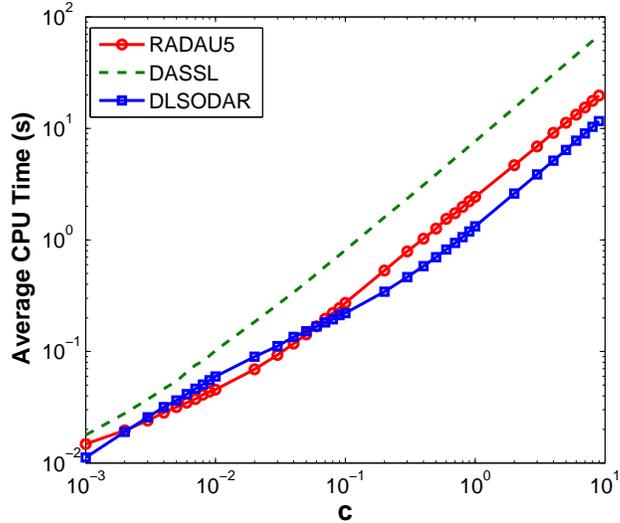}
 \caption{Average CPU times of hybrid method with three ODE solvers on cell cycle model}
 \label{fig:cell_eff}
\end{figure}

To test the influence of slow reactions on the performance of the hybrid method with different ODE solvers, we introduce a scale factor $c$ varying from $10^{-3}$ to $10$, and multiply it to the propensity functions of the six stochastic reactions in the SSA subsystem. With $T=2000$ as the final time, for each set of parameter the system was simulated $1,000$ times using hybrid methods with RADAU5, DASSL, and DLSODAR, respectively.  The mean firing numbers of synthesis and degradation reactions and average CPU time are shown in Figure \ref{fig:cell_acc} and Figure \ref{fig:cell_eff}.

For the comparison of accuracy, we just plot the results for $M_z$ synthesis reaction in Figure \ref{fig:cell_acc}, the propensity function of which contains a hill function. Similar as the previous two examples, hybrid methods with three ODE solvers give similar accuracy results. Figure \ref{fig:cell_eff} shows the average CPU times of hybrid methods with three ODE solvers. DASSL performs the worst compared with the other two ODE solvers. The green dash curve is on the top as $c$ increases. If we focus on the red dot curve (RADAU5) and the blue square curve (DLSODAR),
we note that these two curves have two crossovers, different from the no crossover pattern in Figure \ref{fig:toy_time} and the one crossover pattern in Figure \ref{fig:gene_eff}. When $0.002<c<0.06$, hybrid method with RADAU5 takes less time. But on other scales, hybrid method with DLSODAR has higher efficiency.

\section{Conclusion and discussion}\label{sec:discussion}

The hybrid method which combines both continuous model and stochastic model has received more and more attentions from biochemical system simulation community. It not only keeps the stochastic property of biochemical system but also costs less time to accomplish simulations by dealing with the continuous part with more efficient approaches.  When applying the hybrid method, the continuous model is usually handled by ODE solvers. On the fortran platform, DLSODAR, an ODE solver based on multistep method is widely used. But when integration is interrupted frequently by the stochastic subsystem, the efficiency of mutistep ODE solvers may drop heavily, since the frequent restarting process inhibits ODE solvers to use higher orders, while the one step method ODE solver RADAU5 may have its advantages under this situation.

In this paper, we developed a root finding strategy, which is particularly designed for the hybrid method, and applied them to RADAU5 and DASSL. With the modification, we were able to conduct three groups of numerical experiments to compared the accuracy and efficiency of hybrid methods with three ODE solvers. The results of the first two experiments demonstrate the advantage of RADAU5 when systems contain frequently firing slow reactions. However, in the third experiment on a cell cycle model, DLSODAR has better performances with respect to low and high firing frequency of slow reactions. The reason why RADAU5 no longer holds its advantages for the cell cycle model is intriguing to explore.

RADAU5 is an ODE solver based on the implicit Runge-Kutta method. In order to calculate the state values at the next step, Newton's iteration, the main part of RADAU5, requires accurate Jacobian matrix. Therefore, beside Newton's iteration, the calculation of Jacobian matrix is also an important part, which takes up a large percent of CPU time. Each time when RADAU5 starts, the Jacobian matrix should be calculated immediately. During the integration, if errors are accepted, there is no need to recalculate Jacobian matrix. For the first two models, each set of ODEs involves no more than two species. To obtain the Jacobian matrix is not a hard work. But it is not the same for the cell cycle model. In the cell cycle model, the highly nonlinear terms make the calculation more expensive. With the aid of collected information from running data, we observed that as the restart number grows, RADAU5 spends more and more time on Jacobian matrix, rising from forty percent to sixty percent. However, the multistep method DLSODAR avoids such a process. DLSODAR has a method selection strategy, which automatically switches between Adams methods and BDF methods. When parameter $c$ is beyond a threshold, both RADAU5 and DLSODAR will stop in each step. In such condition, DLSODAR always chooses Adams method to enhance the efficiency. No Jacobian matrix calculation is needed any more. This property largely reduces computation complexity. Therefore, although the average step size of RADAU5 is larger than DLSODAR, the benefit is gradually counteracted by the burden of calculation of Jacobian matrix. That is why the curve referring average CPU time of hybrid method with RADAU5 intercepts the one of hybrid method with DLSODAR twice. When making use of hybrid methods, we should carefully partition systems. If the sum of propensities of slow reaction is large enough, which forces ODE solvers to stop very frequently, DLSODAR may be the first to be considered.

To make the analysis easier in section \ref{sec:inverse}, we made an assumption that $t_z<g_1<g_2<g_3$, since $g_i$ is an approximation to $z(t+c_ih)$ which is an monotonous function. Each time when ODE solvers are interrupted at time $t+h$, we can guarantee $t_z<0$ and $g_3>0$. Of course, the time $t_{event}$ when slow reaction fires must be within time interval $[t, t+h]$. However, in practice, we are not sure what will be the situations for $g_1$ and $g_2$. If values of these four variables have the relationship as we assumed, the inverse interpolation will work very well. If this assumption is not satisfied, inverse interpolation could be dangerous. Under this condition, these four variables no long describe a monotonous function. The result may be far from the true value, for instance we may get $t_{event}>t+h$ or $t_{event}<t$, which is an obvious wrong output. So far we have not observed such big errors appearing in our numerical experiments. But, in order to prevent this bad inverse interpolation from affecting the accuracy heavily, more work need to be done to improve the root funding mechanism.

\vspace{0.2in}
\noindent \textbf{Acknowledgements} This work was partially supported by the
National Science Foundation under awards DMS-1225160, CCF-0953590, and CCF-1526666
and the National Institutes of Health
under award GM078989.


\bibliographystyle{amsplain}
\bibliography{myBib}

\end{document}